\documentclass[a4paper,10pt]{article}
\usepackage{stmaryrd}
\usepackage{amsfonts}
\usepackage{bbm}
\usepackage{amscd}
\usepackage{mathrsfs}
\usepackage{latexsym,amssymb,amsmath,amscd,amscd,amsthm,amsxtra}
\usepackage[dvips]{graphicx}
\usepackage[utf8]{inputenc}
\usepackage[T1]{fontenc}
\usepackage{lmodern}
\usepackage{amssymb}
\usepackage[all]{xy}
\usepackage{nicefrac,mathtools,enumitem}
\usepackage{microtype}
\usepackage{xcolor}
\usepackage{color}
\definecolor{red}{rgb}{1,0,0}
\definecolor{green}{rgb}{0,1,0}
\definecolor{blue}{rgb}{0,0,1}
\definecolor{refkey}{gray}{.625}
\definecolor{labelkey}{gray}{.625}

\newtheorem{thm}{Theorem}[section]
\newtheorem{lem}[thm]{Lemma}
\newtheorem{cor}[thm]{Corollary}
\newtheorem{pro}[thm]{Proposition}
\newtheorem{ex}[thm]{Example}
\newtheorem{rmk}[thm]{Remark}
\newtheorem{defi}[thm]{Definition}

\newcommand {\comment}[1]{{\marginpar{*}\scriptsize\textbf{Comments:} #1}}

\newcommand {\emptycomment}[1]{} 

\newcommand{\be }{\begin{equation}}
\newcommand{\ee }{\end{equation}}

\newcommand{\pf}{\noindent{\bf Proof.}\ }

\newcommand{\huaB}{\mathcal{B}}

\newcommand{\huaA}{\mathcal{A}}
\newcommand{\huaL}{\mathcal{L}}

\newcommand{\huaF}{\mathcal{F}}

\newcommand{\huaP}{\mathcal{P}}

\newcommand{\CWM}{C^{\infty}(M)}

\newcommand{\XM}{\mathcal{X}(M)}

\newcommand{\frkd}{\mathfrak d}
\newcommand{\frke}{\mathfrak e}

\newcommand{\frkg}{\mathfrak g}

\newcommand{\frkD}{\mathfrak D}

\newcommand{\frkP}{\mathfrak P}

\def\qed{\hfill ~\vrule height6pt width6pt depth0pt}

\newcommand{\half}{\frac{1}{2}}

\newcommand{\br}[1]{   [ \cdot,    \cdot  ]   }

\newcommand{\dev}{\mathfrak{D}}

\newcommand{\id}{\rm{id}}

\newcommand{\g}{\mathfrak g}

\newcommand{\dM}{\mathrm{d}}

\newcommand{\Hom}{\mathrm{Hom}}

\newcommand{\gl}{\mathfrak {gl}}

\newcommand{\End}{\mathrm{End}}
\newcommand{\ad}{\mathrm{ad}}

\newcommand{\Img}{\mathrm{Im}}

\begin{document}
\title{Pre-symplectic   algebroids and their applications\thanks
 {
This research  is supported by
NSF of China (11471139, 11271202, 11221091, 11425104), SRFDP
(20120031110022) and NSF of Jilin Province (20140520054JH).
 }}
\author{Jiefeng Liu$^1$, Yunhe Sheng$^{1,2}$ and Chengming Bai$^3$\\
$^1$Department of Mathematics, Xinyang Normal University,\\ \vspace{2mm} Xinyang 464000, Henan, China\\
 $^2$Department of Mathematics, Jilin University,
 \\\vspace{2mm}Changchun 130012, Jilin, China\\
$^3$Chern Institute of Mathematics and LPMC, Nankai University,\\
Tianjin 300071, China \\ Email:  jfliu12@126.com; shengyh@jlu.edu.cn; baicm@nankai.edu.cn \\
}

\date{}

\maketitle
\begin{abstract}
In this paper, we introduce the notion of a pre-symplectic algebroid, and show that there is a one-to-one correspondence between pre-symplectic algebroids and symplectic Lie algebroids. This result is the geometric generalization of the relation between left-symmetric algebras and symplectic (Frobenius) Lie algebras. Although pre-symplectic algebroids are not left-symmetric algebroids, they still can be viewed as the underlying structures of symplectic Lie algebroids. 
Then we study exact pre-symplectic algebroids and show that they are classified by the third cohomology group of a left-symmetric algebroid. Finally, we study para-complex pre-symplectic algebroids. Associated to a para-complex pre-symplectic algebroid, there is a pseudo-Riemannian Lie algebroid. The multiplication in a para-complex pre-symplectic algebroid characterizes the restriction to the Lagrangian subalgebroids of the Levi-Civita connection in the corresponding pseudo-Riemannian Lie algebroid.

\end{abstract}

\section{Introduction}

  Left-symmetric algebras (or pre-Lie algebras) arose from the
study of convex homogeneous cones \cite{Vinb}, affine manifolds
and affine structures on Lie groups \cite{Koszul1}, deformation
and cohomology theory of associative algebras \cite{G} and then
appear in many fields in mathematics and mathematical physics. See
the survey article \cite{Pre-lie algebra in geometry} and the
references therein. In particular, there are
 close relations between left-symmetric algebras and
certain important left-invariant structures on Lie groups like
aforementioned affine, symplectic, K\"ahler, and metric structures
\cite{Kim,Lichnerowicz,Med,Milnor}.

 On the other hand, the notion of a Lie algebroid was introduced
by Pradines in 1967, which is a generalization of Lie algebras and
tangent bundles. See \cite{General theory of Lie groupoid and Lie
algebroid} for general theory about Lie algebroids. They play
important roles in various parts of mathematics.
 In particular, the notion of a symplectic Lie algebroid was introduced in \cite{NestTsygan}, and the geometric construction of deformation quantization of a  symplectic Lie algebroid was studied there. There are important applications of symplectic Lie algebroids in  geometric mechanics \cite{symplectic lie algebroid,reduction,reduction1}. A symplectic Lie algebroid is a Lie algebroid equipped with a nondegenerate 2-cocycle. It can be viewed as the geometric generalization of a symplectic (Frobenius) Lie algebra, which is a Lie algebra $\g$ equipped with a nondegenerate 2-cocycle $\omega\in\wedge^2\g^*$.  
There is a one-to-one correspondence between symplectic (Frobenius) Lie algebras and quadratic left-symmetric algebras, which are left-symmetric algebras together with nondegenerate invariant skew-symmetric bilinear forms \cite{symplectic Lie algebras}. The relation is given as follows:
$$
\xymatrix{
 \mbox{symplectic~ (Frobenius)~ Lie~ algebra}\ar@<1ex>[rr]^{x\cdot y={\omega^\sharp}^{-1}\ad^*_x\omega^\sharp(y)}
                && \mbox{quadratic ~left-symmetric  ~algebra.} \ar[ll]^{\mbox{commutator}}  }
$$

In  \cite{LiuShengBaiChen}, we introduced  the notion of a
left-symmetric algebroid\footnote{Recently, we realize that this structure was already studied by Nguiffo Boyom in \cite{Boyom1} under the name of  Koszul-Vinberg algebroids. See also \cite{Boyom2} for more details about the relation with Poisson manifolds.}, which is a geometric generalization of a
left-symmetric algebra.\emptycomment{, such that we have the following
commutative diagram:
$$
\xymatrix{
 \mbox{Lie~ algebras} \ar[d]_{\mbox{geometrization}}
                && \mbox{left-symmetric  ~algebras} \ar[d]^{\mbox{geometrization}}\ar[ll]_{\mbox{commutator}\qquad}  \\
\mbox{Lie~ algebroids}
                && \mbox{left-symmetric  ~algebroids.}    \ar[ll]_{\mbox{commutator}\qquad}    }
$$$$
            \footnotesize{ \mbox{ Diagram~ I}}
$$
}
We developed the basic theory of left-symmetric algebroids and gave applications in Poisson algebras, symplectic Lie algebroids and complex Lie algebroids.

In this paper, first we introduce the notion of a pre-symplectic algebroid
which is a geometric generalization of a quadratic left-symmetric
algebra. We would like to point out that one of the nontriviality is
that a pre-symplectic algebroid is not a left-symmetric algebroid.
More importantly, there is a one-to-one correspondence between
symplectic Lie algebroids and pre-symplectic algebroids (Theorem
\ref{symp 3} and Theorem \ref{thm:Symplectic-LWX}). This result
generalizes the relation between symplectic~ (Frobenius)~ Lie~
algebras and quadratic  left-symmetric   algebras. In terms of
 commutative diagram, we have
$$
\xymatrix{
 \mbox{symplectic~ (Frobenius)~ Lie~ algebra}\ar@<1ex>[rr]^{x\cdot y={\omega^\sharp}^{-1}\ad^*_x\omega^\sharp(y)} \ar[d]_{\mbox{geometrization}}
                && \mbox{quadratic ~left-symmetric  ~algebra} \ar[d]^{\mbox{geometrization}}\ar[ll]^{\mbox{commutator}}  \\
\mbox{symplectic~Lie~ algebroid}\ar@<1ex>[rr]^{\footnotesize{e_1\star e_2={\omega^\sharp}^{-1}(\huaL_{e_1}\omega^\sharp(e_2)+\half \dM(\omega(e_1,e_2)))}}
                && \mbox{pre-symplectic algebroid.}    \ar[ll]^{\mbox{commutator}} }
$$
$$
 \footnotesize{ \mbox{ Diagram~ I}}
$$
A pre-symplectic algebroid is defined to be a quadruple
$(E,\star,\rho,(\cdot,\cdot)_-)$, where $E$ is a vector bundle
over $M$, $\star$ is a multiplication on the section space
$\Gamma(E)$, $\rho:E\longrightarrow TM$ is a bundle map and
$(\cdot,\cdot)_-$ is a nondegenerate skew-symmetric bilinear form
such that some compatibility conditions are satisfied.
Even though a pre-symplectic  algebroid is not a left-symmetric
algebroid, its Dirac structures are all left-symmetric algebroids.
The notion of a Dirac structure was invented to treat in the same
framework of Poisson structures and symplectic structures
\cite{Courant}. There is a one-to-one correspondence between Dirac
structures in a pre-symplectic algebroid and Langrangian
subalgebroids in the corresponding symplectic Lie algebroid.
Moreover, associated to any left-symmetric algebroid $A$, there is
a canonical pre-symplectic algebroid structure on $A\oplus A^*$,
which is called the pseudo-semidirect product of $A$ and $A^*$.

Then we study exact pre-symplectic algebroids $(E,\star,\rho,(\cdot,\cdot)_-)$ with respect to the left-symmetric algebroid $T_\nabla M$ associated to a flat manifold $(M,\nabla)$. By choosing an isotropic splitting, $E$ is isomorphic to $TM\oplus T^*M$ as vector bundles. The pre-symplectic algebroid structure on $TM\oplus T^*M$ is a twist of the pseudo-semidirect product by a 3-cocycle. Consequently, we classify exact pre-symplectic algebroids by the third cohomology group of the left-symmetric algebroid $T_\nabla M$. This result is parallel to \v Severa's result on the classification of exact Courant algebroids \cite{severa3form}.

Finally we study para-complex pre-symplectic algebroids, which are
the underlying structures of para-K\"{a}hler Lie algebroids.
Para-K\"{a}hler Lie algebroids were introduced by Leichtnam, Tang
and Weinstein in \cite{Para-kahler} to describe the deformation
quantization near a strictly pseudoconvex boundary. See
\cite{Benayadi,Cruc2,Kan} for more details for para-K\"{a}hler Lie algebras and para-K\"{a}hler manifolds.
In the recent work \cite{Vai06}, Vaisman studied generalized
para-K\"{a}hler manifolds. In this paper, we introduce the notion
of a para-complex pre-symplectic algebroid. There is a natural
example of para-complex pre-symplectic algebroid constructed from a
pseudo-semidirect product. We prove that there is a one-to-one
correspondence between para-K\"{a}hler Lie algebroids and
para-complex pre-symplectic algebroids (Proposition
\ref{pro:PC-PK}). Given a para-K\"{a}hler Lie algebroid, there is
an associated  pseudo-Riemannian Lie algebroid  satisfying some
compatibility conditions. We find that the multiplication  in a
para-complex pre-symplectic algebroid
$(E,\star,\rho,(\cdot,\cdot)_-,P)$ coincides with the Levi-Civita
$E$-connection in the corresponding pseudo-Riemannian Lie
algebroid  $(E,[\cdot,\cdot]_E,\rho,g)$ restricting to the
Lagrangian subalgebroids (Proposition \ref{pro:mul-Levi}). This
provides potential applications of pre-symplectic algebroids in
these areas.

The paper is organized as follows. In Section $2$, we give a review on left-symmetric algebras,  Lie algebroids and left-symmetric algebroids. In Section $3$,  we introduce the notion of pre-symplectic algebroids and give their relation with symplectic Lie algebroids. In  Section $4$, we study exact pre-symplectic algebroids and give their classification. In Section $5$, we study para-complex pre-symplectic algebroids and give their applications.

Throughout this paper, all the algebras and vector spaces
are assumed to be finite-dimensional, although many results still
hold in the infinite-dimensional cases. All the vector bundles are over the same manifold
$M$. For two vector bundles $A$ and $B$, a bundle map from $A$ to
$B$ is a base-preserving map and $\CWM$-linear.

\vspace{2mm}
 \noindent {\bf Acknowledgement:} We give our warmest thanks to Zhangju Liu and Jianghua Lu for very useful comments and discussions.

\section{Preliminaries}
\subsection*{$\bullet$ Left-symmetric algebras}
\begin{defi}
A {\bf left-symmetric algebra} is a pair $(\frkg,\cdot_\frkg)$, where $\g$ is a vector space, and  $\cdot_\frkg:\g\otimes \g\longrightarrow\g$ is a bilinear multiplication
satisfying that for all $x,y,z\in \g$, the associator
\begin{equation}\label{eq:associator}
(x,y,z)\triangleq x\cdot_\frkg(y\cdot_\frkg z)-(x\cdot_\frkg y)\cdot_\frkg z
\end{equation} is symmetric in $x,y$,
i.e.
$$(x,y,z)=(y,x,z),\;\;{\rm or}\;\;{\rm
equivalently,}\;\;x\cdot_\frkg(y\cdot_\frkg z)-(x\cdot_\frkg y)\cdot_\frkg z=y\cdot_\frkg(x\cdot_\frkg z)-(y\cdot_\frkg x)\cdot_\frkg
z.$$
\end{defi}
For all $x\in \g$, let $L_x$
 and $R_x$ denote the left and right multiplication operators respectively, i.e. $$L_x(y)=x\cdot_\frkg
y,\quad R_xy=y\cdot_\g x,\quad \forall ~y\in \g.$$
The following conclusion is known
\cite{Pre-lie algebra in geometry}:
\begin{lem} \label{lem:sub-ad} Let $(\g,\cdot_\g)$ be a left-symmetric algebra. The commutator
$ [x,y]_\g=x\cdot_\g y-y\cdot_\g x$ defines a Lie algebra $\g^c$,
which is called the {\bf sub-adjacent Lie algebra} of $\g$ and $\g$ is
called a {\bf compatible left-symmetric algebra}  on the Lie
algebra $\g^c$. Furthermore, $L:\g^c\rightarrow
\gl(\g)$ gives a representation of the Lie
algebra $\g^c$.
\end{lem}

A nondegenerate  skew-symmetric bilinear form
$(\cdot,\cdot)_-$ on a left-symmetric algebra $(\frkd,\cdot_\frkd)$
is called {\bf invariant} if
\begin{equation}
(a\cdot_\frkd b,c)_-+(b,[a,c]_\frkd)_-=0.
\end{equation}
A {\bf quadratic left-symmetric algebra} is a left-symmetric algebra equipped with a nondegenerate skew-symmetric invariant bilinear form.

\subsection*{$\bullet$ Lie algebroids}
\begin{defi}
A {\bf Lie algebroid} structure on a vector bundle $\huaA\longrightarrow M$ is
a pair that consists of a Lie algebra structure $[\cdot,\cdot]_\huaA$ on
the section space $\Gamma(\huaA)$ and a  bundle map
$a_\huaA:\huaA\longrightarrow TM$, called the anchor, such that the
following relation is satisfied:
$$~[x,fy]_\huaA=f[x,y]_\huaA+a_\huaA(x)(f)y,\quad \forall~f\in
\CWM.$$
\end{defi}

For a vector bundle $E\longrightarrow M$, let
$\dev(E)$ be the gauge Lie algebroid of the
 frame bundle
 $\huaF(E)$, which is also called the covariant differential operator bundle of $E$.

Let $(\huaA,[\cdot,\cdot]_\huaA,a_\huaA)$ and $(\huaB,[\cdot,\cdot]_\huaB,a_\huaB)$ be two Lie
algebroids (with the same base), a {\bf base-preserving morphism}
from $\huaA$ to $\huaB$ is a bundle map $\sigma:\huaA\longrightarrow \huaB$ such
that
\begin{eqnarray*}
  a_\huaB\circ\sigma=a_\huaA,\quad
  \sigma[x,y]_\huaA=[\sigma(x),\sigma(y)]_\huaB.
\end{eqnarray*}

A {\bf representation} of a Lie algebroid $\huaA$
 on a vector bundle $E$ is a base-preserving morphism $\rho$ form $\huaA$ to the Lie algebroid $\dev(E)$.
Denote a representation by $(E;\rho).$
The {\bf dual representation} of the Lie algebroid $\huaA$ on $E^*$ is the bundle map $\rho^*:\huaA\longrightarrow \dev(E^*)$ given by
\begin{equation}\label{eq:dual}
\langle \rho^*(x)(\xi),y\rangle=a_\huaA(x)\langle \xi,y\rangle-\langle \xi,\rho(x)(y)\rangle,\quad \forall~x,y\in \Gamma(\huaA),~\xi\in\Gamma(E^*).
\end{equation}

For all $x\in \Gamma(\huaA)$, the {\bf Lie derivation} $\huaL_x:\Gamma(\huaA^*)\longrightarrow\Gamma(\huaA^*)$ of the Lie algebroid $\huaA$  is given by
\begin{eqnarray}\label{Lie der1}
\langle\huaL_x\xi,y\rangle=a_\huaA(x)\langle \xi,y\rangle-\langle \xi,[x,y]_\huaA\rangle,\quad \forall y\in\Gamma(\huaA),\xi\in\Gamma(\huaA^*).
\end{eqnarray}

A Lie algebroid $(\huaA,[\cdot,\cdot]_\huaA,a_\huaA)$ naturally
represents on the trivial line bundle $E=M\times \mathbb R$ via
the anchor map $a_\huaA:\huaA\longrightarrow TM$. The
corresponding coboundary operator
$\dM:\Gamma(\wedge^k\huaA^*)\longrightarrow
\Gamma(\wedge^{k+1}\huaA^*)$ is given by
\begin{eqnarray*}
  \dM\varpi(x_1,\cdots,x_{k+1})&=&\sum_{i=1}^{k+1}(-1)^{i+1}a_\huaA(x_i)\varpi(x_1\cdots,\widehat{x_i},\cdots,x_{k+1})\\
  &&+\sum_{i<j}(-1)^{i+j}\varpi([x_i,x_j]_\huaA,x_1\cdots,\widehat{x_i},\cdots,\widehat{x_j},\cdots,x_{k+1}).
\end{eqnarray*}
In particular, a $2$-form $\varpi\in\Gamma(\wedge^2\huaA^*)$ is a {\bf 2-cocycle} if $\dM \varpi=0$, i.e.
\begin{equation}
  a_\huaA(x)\varpi(y,z)- a_\huaA(y)\varpi(x,z)+ a_\huaA(z)\varpi(x,y)-\varpi([x,y]_\huaA,z)+\varpi([x,z]_\huaA,y)-\varpi([y,z]_\huaA,x)=0.
\end{equation}

A {\bf symplectic Lie algebroid} is a Lie algebroid together with a nondegenerate closed $2$-form. A subalgebroid of a symplectic Lie algebroid $(\huaA,[\cdot,\cdot]_\huaA,a_\huaA,\varpi)$ is called {\bf Lagrangian} if it is maximal isotropic with respect to the skew-symmetric bilinear form $\varpi$.

\subsection*{$\bullet$ Left-symmetric algebroids}
\begin{defi}\label{defi:left-symmetric algebroid}
A {\bf left-symmetric algebroid} structure on a vector bundle
$A\longrightarrow M$ is a pair that consists of a left-symmetric
algebra structure $\cdot_A$ on the section space $\Gamma(A)$ and a
vector bundle morphism $a_A:A\longrightarrow TM$, called the anchor,
such that for all $f\in\CWM$ and $x,y\in\Gamma(A)$, the following
conditions are satisfied:
\begin{itemize}
\item[\rm(i)]$~x\cdot_A(fy)=f(x\cdot_A y)+a_A(x)(f)y,$
\item[\rm(ii)] $(fx)\cdot_A y=f(x\cdot_A y).$
\end{itemize}
\end{defi}

We usually denote a left-symmetric algebroid by $(A,\cdot_A, a_A)$.

A left-symmetric algebroid $(A,\cdot_A,a_A)$ is called {\bf transitive (regular)} if  $a_A$ is
 surjective (of constant rank). Any left-symmetric  algebra is a left-symmetric algebroid over a point.

 \begin{ex}\label{ex:main}{\rm Let $M$
be a differential manifold with a flat torsion free connection
$\nabla$. Then $(TM,\nabla,\id)$ is a left-symmetric algebroid whose sub-adjacent Lie algebroid is
exactly the tangent Lie algebroid. We denote this left-symmetric algebroid by $T_\nabla M$, which will be frequently used below.
}
\end{ex}

For any $x\in\Gamma(A)$, we define
$L_x:\Gamma(A)\longrightarrow\Gamma(A)$ by $L_x(y)=x\cdot_A y$ and $R_x:\Gamma(A)\longrightarrow\Gamma(A)$ by $R_x(y)=y\cdot_A x$ .
Condition (i) in the above definition means that $L_x\in \frkD(A)$.
Condition (ii) means that the map $x\longmapsto L_x$ is
$C^\infty(M)$-linear. Thus, $L:A\longrightarrow \frkD(A)$ is a bundle
map. Furthermore, condition (ii) also implies that  $R_x:\Gamma(A)\longrightarrow\Gamma(A)$ is $C^\infty(M)$-linear. Define $R^*_x:\Gamma(A^*)\longrightarrow\Gamma(A^*)$ by
\begin{equation}\label{eq:dualR}
\langle R_x^*\xi,y\rangle=-\langle \xi,R_x y\rangle,\quad \forall x,y\in\Gamma(A), \xi\in\Gamma(A^*).
\end{equation}
\begin{pro}{\rm \cite{LiuShengBaiChen}}\label{thm:sub-adjacent}
  Let $(A,\cdot_A, a_A)$ be a left-symmetric algebroid. Define  a skew-symmetric bilinear bracket operation $[\cdot,\cdot]_A$ on $\Gamma(A)$ by
  $$
  [x,y]_A=x\cdot_A y-y\cdot_A x,\quad \forall ~x,y\in\Gamma(A).
  $$
Then, $(A,[\cdot,\cdot]_A,a_A)$ is a Lie algebroid, and denoted by
$A^c$, called the {\bf sub-adjacent Lie algebroid} of
 $(A,\cdot_A,a_A)$. Furthermore, $L:A\longrightarrow \frkD(A)$  gives a
  representation of the Lie algebroid  $A^c$.
\end{pro}

\begin{thm}{\rm \cite{LiuShengBaiChen}}\label{thm:symLie}
   Let $(A,\cdot_A, a_A)$ be a left-symmetric algebroid. Then $(A^c\ltimes_{L^*}A^*,[\cdot,\cdot]_S,\rho,\omega)$ is a symplectic Lie algebroid, where $A^c\ltimes_{L^*}A^*$ is the semidirect product of $A^c$ and $A^*$ in which $L^*$ is the dual representation\footnote{Note that the definitions of $L^*$ (see \eqref{eq:dual}) and the above defined $R^*$ (see \eqref{eq:dualR}) are not the same.} of $L$. More precisely, the Lie bracket $[\cdot,\cdot]_S$ and the anchor $\rho$ are given by
   $$
   [x+\xi,y+\eta]_S=[x,y]_A+L^*_x\eta-L^*_y\xi,
   $$
 and $\rho(x+\xi)=a_A(x)$  respectively.
   Furthermore, the symplectic form $\omega$ is given by
   \begin{equation}\label{eq:defiomega}
     \omega(x+\xi,y+\eta)=\langle\xi,y\rangle-\langle\eta,x\rangle,\quad \forall x,y\in\Gamma(A),~\xi,\eta\in\Gamma(A^*).
   \end{equation}
\end{thm}




Let $(A,\cdot_A,a_A)$ be a left-symmetric algebroid and $E$  a vector
bundle. A {\bf representation} of $A$ on $E$ consists of a pair
$(\rho,\mu)$, where $\rho:A\longrightarrow \frkD(E)$ is a representation
of $A^c$ on $E $ and $\mu:A\longrightarrow \End(E)$ is a bundle
map, such that for all $x,y\in \Gamma(A),\ e\in\Gamma(E)$, we have
\begin{eqnarray}\label{representation condition 2}
 \rho(x)\mu(y)e-\mu(y)\rho(x)e=\mu(x\cdot_A y)e-\mu(y)\mu(x)e.
\end{eqnarray}
Denote a representation by $(E;\rho,\mu)$.

\emptycomment{let us recall the cohomology complex for a left-symmetric algebroid $(A,\cdot_A,a_A)$ with a representation $(E;\rho,\mu)$ briefly. Denote the set of $(n+1)$-cochains by
$$C^{n+1}(A,E)=\Gamma(\Hom(\wedge^{n}A\otimes A,E)),\
n\geq 0.$$  For all $\varphi\in C^{n}(A,E)$, the coboundary operator $\delta:C^{n}(A,E)\longrightarrow C^{n+1}(A,E)$ is given by
 \begin{eqnarray*}
\delta\varphi(x_1,\cdots,x_{n+1})&=&\sum_{i=1}^{n}(-1)^{i+1}\rho(x_i)\omega(x_1,\cdots,\hat{x_i},\cdots,x_{n+1})\\
 &&+\sum_{i=1}^{n}(-1)^{i+1}\mu(x_{n+1})\omega(x_1,\cdots,\hat{x_i},\cdots,x_n,x_i)\\
 &&-\sum_{i=1}^{n}(-1)^{i+1}\varphi(x_1,\cdots,\hat{x_i},\cdots,x_n,x_i\cdot_A x_{n+1})\\
 &&+\sum_{1\leq i<j\leq n}(-1)^{i+j}\varphi([x_i,x_j]_A,x_1,\cdots,\hat{x_i},\cdots,\hat{x_j},\cdots,x_{n+1}),
\end{eqnarray*}
for all $x_i\in \Gamma(A),i=1,\cdots,n+1$.
}

At the end of this section,  let us recall the cohomology complex with the coefficients in the trivial representation, i.e. $\rho=a_A$ and $\mu=0$. See \cite{cohomology of pre-Lie,LiuShengBaiChen} for general theory of cohomologies of right-symmetric algebras and left-symmetric algebroids respectively. The set of $(n+1)$-cochains  is given by
$$C^{n+1}(A)=\Gamma(\wedge^{n}A^*\otimes A^*),\
n\geq 0.$$  For all $\varphi\in C^{n}(A)$ and $x_i\in
\Gamma(A),~i=1,\cdots,n+1$, the corresponding coboundary operator
$\delta$ is given by
 \begin{eqnarray}\label{LSCA cohomology}
\nonumber\delta\varphi(x_1, \cdots,x_{n+1})
 &=&\sum_{i=1}^{n}(-1)^{i+1}a_A(x_i)\varphi(x_1, \cdots,\hat{x_i},\cdots,x_{n+1})\nonumber\\
 &&-\sum_{i=1}^{n}(-1)^{i+1}\varphi(x_1, \cdots,\hat{x_i},\cdots,x_n,x_i\cdot_A x_{n+1})\nonumber\\
 &&+\sum_{1\leq i<j\leq n}(-1)^{i+j}\varphi([x_i,x_j]_A,x_1,\cdots,\hat{x_i},\cdots,\hat{x_j},\cdots,x_{n+1}).
\end{eqnarray}


 \section{Pre-symplectic  algebroids and symplectic Lie algebroids}

 In this section, we introduce the notion of a pre-symplectic algebroid. We show that there is a one-to-one correspondence between pre-symplectic algebroids and symplectic Lie algebroids.
  Thus, pre-symplectic algebroids can be viewed as the underlying structures of symplectic Lie algebroids.
 Furthermore, there is a natural
construction of pre-symplectic algebroids from  left-symmetric
algebroids (the so-called pseudo-semidirect products).

 \begin{defi}\label{MTLA}
 A {\bf pre-symplectic algebroid} is a vector bundle $E\rightarrow M$ equipped with a nondegenerate skew-symmetric bilinear form $(\cdot,\cdot)_-$, a multiplication $\star:\Gamma(E)\times\Gamma(E)\longrightarrow\Gamma(E) $, and a bundle map $\rho:E\rightarrow TM$, such that  for all $e, e_1,e_2,e_3\in\Gamma(E),~~f\in C^\infty(M)$, the following conditions are satisfied:
\begin{itemize}
\item[$\rm(i)$] $(e_1,e_2,e_3)-(e_2,e_1,e_3)=\frac{1}{6}DT(e_1,e_2,e_3)$\label{associator};
\item[$\rm(ii)$]$\rho(e_1)(e_2,e_3)_-=(e_1{\star}e_2-\frac{1}{2}D(e_1,e_2)_-,e_3)_-+(e_2,[e_1,e_3]_E)_-,$\label{invariant}
   \end{itemize}
   where $(e_1,e_2,e_3)$ is the associator for the multiplication $\star$ given by \eqref{eq:associator}, $T:\Gamma(E)\times\Gamma(E)\times\Gamma(E) \longrightarrow \CWM$ is defined by
   \begin{equation}\label{T-equation}
   T(e_1,e_2,e_3)= (e_1\star e_2, e_3)_-+(e_1,e_2\star e_3)_- - (e_2\star e_1, e_3)_--(e_2,e_1\star e_3)_-,
   \end{equation}
 $D:C^\infty(M)\longrightarrow\Gamma(E)$ is defined by
  \begin{equation}\label{eq:defD}
 (Df,e)_-=\rho(e)(f),
 \end{equation}
 and the bracket $[\cdot,\cdot]_E:\wedge^2\Gamma(E)\longrightarrow \Gamma(E)$ is defined by
 \begin{equation}\label{eq:defbraE}
 [e_1,e_2]_E=e_1\star e_2-e_2\star e_1.
 \end{equation}
\end{defi}

We denote a pre-symplectic algebroid  by  $(E,\star,\rho,(\cdot,\cdot)_-)$.
\begin{rmk}
  When $M$ is a point,  we recover the notion of a quadratic left-symmetric algebra. It is well-known
that there is a one-to-one correspondence between quadratic left-symmetric
algebras and symplectic (Frobenius) Lie algebras. Note that a symplectic Lie algebroid is a geometric generalization of a symplectic
Lie algebra. In the following we will show that a pre-symplectic algebroid can give rise to
a symplectic Lie algebroid. Conversely, one can also obtain a pre-symplectic algebroid
from a symplectic Lie algebroid. This justifies that a pre-symplectic algebroid is a geometric generalization of a  quadratic left-symmetric algebra as given in Diagram I.
\end{rmk}

\begin{rmk}
 A pre-symplectic algebroid, as a geometric generalization of a
quadratic left-symmetric algebra, is not simply a left-symmetric
algebroid with a nondegenerate skew-symmetric bilinear form
satisfying certain compatibility conditions. This is quite interesting.
\end{rmk}

\begin{rmk}
 In this remark, we point out that the algebraic structure underlying a pre-symplectic algebroid is a $2$-term pre-Lie$_{\infty}$-algebra, or equivalently a pre-Lie $2$-algebra \cite{sheng}. See \cite{chapoton} for more details about pre-Lie$_{\infty}$-algebras. More precisely, let $(E,\star,\rho,(\cdot,\cdot)_-)$ be a pre-symplectic algebroid. Consider the graded vector space $\frke=\frke_0\oplus \frke_1$, where $ \frke_0=\Gamma(E)$ and $\frke_1= \Img(D)$,
 then $(\frke,l_1,l_2,l_3)$ is a $2$-term pre-Lie$_{\infty}$-algebra, where $l_i$ are given by the following formulas:
\begin{eqnarray*}
l_1(x)&=&x,\quad \forall x\in\frke_1,\\
l_2(e_1\otimes e_2)&=&e_1\star e_2,\quad \forall e_1,e_2\in\frke_0,\\
l_2(e\otimes x)&=&e\star x,\quad \forall e\in\frke_0,x\in\frke_1,\\
l_2(x\otimes e)&=&x\star e,\quad \forall e\in\frke_0,x\in\frke_1,\\
l_3(e_1\wedge e_2\otimes e_3)&=&\frac{1}{6} DT(e_1,e_2,e_3),\quad \forall e_1,e_2,e_3\in\frke_0,
\end{eqnarray*}
where $T$ is given by \eqref{T-equation}.
 \end{rmk}

It is obvious that the operator $D$ defined by \eqref{eq:defD} satisfies the Leibniz rule:
 \begin{equation}\label{Leib 3}
 D(fg)=fD(g)+gD(f),\quad ~~\forall f,g\in C^\infty(M).
   \end{equation}

\begin{pro}\label{Condii iii}
Let $(E,\star,\rho,(\cdot,\cdot)_-)$ be a pre-symplectic algebroid. For all $e,e_1,e_2\in\Gamma(E),~~f\in\CWM$, we have
\begin{eqnarray}
e_1\star fe_2&=&f(e_1\star e_2)+\rho(e_1)(f)e_2+\frac{1}{2}(e_1,e_2)_-Df,\label{prop 1}\\
\label{eq:LWXanchor}~[e_1,fe_2]_E&=&f[e_1,e_2]_E+\rho(e_1)(f)e_2,\\
(fe_1)\star e_2&=&f(e_1\star e_2)-\frac{1}{2}(e_1,e_2)_-Df.\label{prop 2}
\end{eqnarray}
\end{pro}
\pf On one hand, by condition (ii) in Definition \ref{MTLA}, we have
    \begin{eqnarray}
     \label{eq:t1} \rho(e_1)(fe_2,e_3)_-&=&(e_1\star fe_2-\frac{1}{2}D(e_1,fe_2)_-,e_3)_-+(fe_2,[e_1,e_3]_E)_-,\\
      \label{eq:t2}\rho(e_1)(e_2,fe_3)_-&=&f(e_1\star e_2-\frac{1}{2}D(e_1,e_2)_-,e_3)_-+(e_2,[e_1,fe_3]_E)_-.
    \end{eqnarray}
On the other hand, the Leibniz rule gives
\begin{eqnarray}
\nonumber\rho(e_1)(fe_2,e_3)_-&=&\rho(e_1)(e_2,fe_3)_-=f\rho(e_1)(e_2,e_3)_-+(e_2,e_3)_-\rho(e_1)(f)\\
\label{eq:t3}&=&f(e_1{\star}e_2-\frac{1}{2}D(e_1,e_2)_-,e_3)_-+f(e_2,[e_1,e_3]_E)_-+(e_2,e_3)_-\rho(e_1)(f).
\end{eqnarray}
By \eqref{eq:t1} and \eqref{eq:t3}, we have
$$(e_1\star fe_2-\frac{1}{2}D(e_1,fe_2)_-,e_3)_-=f(e_1{\star}e_2-\frac{1}{2}D(e_1,e_2)_-,e_3)_-+(e_2,e_3)_-\rho(e_1)(f).$$
Since $(\cdot,\cdot)_-$ is nondegenerate, we have
$$e_1\star fe_2=f(e_1\star e_2)+\rho(e_1)(f)e_2+\frac{1}{2}(e_1,e_2)_-Df,$$
which implies that $(\ref{prop 1})$ holds.

By \eqref{eq:t2} and \eqref{eq:t3}, we can get
 \begin{eqnarray*}
 [e_1,fe_3]_E=f[e_1,e_3]_E+\rho(e_1)(f)e_3,
 \end{eqnarray*}
which implies that \eqref{eq:LWXanchor} holds.

Finally, by \eqref{eq:defbraE}, $(\ref{prop 1})$ and \eqref{eq:LWXanchor}, we can deduce that  $(\ref{prop 2})$ holds. \qed

\emptycomment{
In following, we show the equation $(\ref{prop 3})$ holds.
If $e=0$ the identity is trivial. We assume $y\neq0$.
For all $e,e_1\in\Gamma(E),~~f\in C^\infty(M)$, (ii) gives
$$\rho(fe_1)(e,e)_-=((fe_1)\star e-\frac{1}{2}D(fe_1,e)_-,e)_-+(e,[fe_1,e]_E).$$
Since the $(\cdot,\cdot)_-$ is skew-symmetric and the Leibniz rule of $D$, we have
\begin{eqnarray*}
((e_1,e)_-Df,e)_-+(\rho(e)(f)e_1,e)_-=0
\end{eqnarray*}
Since the $(\cdot,\cdot)_-$ is nondegenerate, we have
$$(Df,e)_-=\rho(e)(f).\qed$$
}

\begin{lem}\label{Jacob Id}
Let $(E,\star,\rho,(\cdot,\cdot)_-)$ be a pre-symplectic algebroid. Then we have
\begin{equation}
T(e_1,e_2,e_3)+c.p.=0,\quad \forall e_1,e_2,e_3\in\Gamma(E),
\end{equation}
where $T$ is given by \eqref{T-equation}.
\end{lem}
\pf It follows from straightforward calculation.\qed

\begin{thm}\label{symp 3}
Let $(E,\star,\rho,(\cdot,\cdot)_-)$ be a pre-symplectic algebroid. Then
 $(E,[\cdot,\cdot]_E,\rho,\omega=(\cdot,\cdot)_-)$ is a symplectic Lie algebroid.
\end{thm}
\pf First, we have
$$[e_1,[e_2,e_3]_E]_E+c.p.=\frac{1}{6}DT(e_1,e_2,e_3)+c.p..$$
By Lemma $\ref{Jacob Id}$, we have
$$[e_1,[e_2,e_3]_E]_E+c.p.=0.$$
By \eqref{eq:LWXanchor}, we have
$$[e_1,fe_2]_E=f[e_1,e_2]_E+\rho(e_1)(f)e_2.$$
Therefore, $(E,[\cdot,\cdot]_E,\rho)$ is a Lie algebroid.

 By condition (ii) in Definition \ref{MTLA}, we have
\begin{eqnarray*}
\rho(e_1)(e_2,e_3)_-&=&(e_1{\star}e_2-\frac{1}{2}D(e_1,e_2)_-,e_3)_-+(e_2,[e_1,e_3]_E)_-;\\
\rho(e_2)(e_1,e_3)_-&=&(e_2{\star}e_1-\frac{1}{2}D(e_2,e_1)_-,e_3)_-+(e_1,[e_2,e_3]_E)_-.
\end{eqnarray*}
Therefore, we have
\begin{eqnarray*}
&&\rho(e_1)(e_2,e_3)_--\rho(e_2)(e_1,e_3)_-\\
&=&(e_1{\star}e_2- e_2{\star}e_1-D(e_1,e_2)_-,e_3)_-+(e_2,[e_1,e_3]_E)_--(e_1,[e_2,e_3]_E)_-\\
&=&([e_1,e_2]_E,e_3)_--([e_1,e_3]_E,e_2)_-+([e_2,e_3]_E,e_1)_--\rho(e_3)(e_1,e_2)_-,
\end{eqnarray*}
which implies that
$$\dM\omega(e_1,e_2,e_3)=0.$$
Thus,  $(E,[\cdot,\cdot]_E,\rho,\omega)$  is a symplectic Lie algebroid.  \qed

\begin{cor}\label{eq:preLie2a1}
Let $(E,\star,\rho,(\cdot,\cdot)_-)$ be a pre-symplectic algebroid. For all $e\in\Gamma(E)$, $f\in\CWM$, we have
\begin{equation}\label{pro1 MLA}
e\star Df=\frac{1}{2}D(Df,e)_-.
\end{equation}
\end{cor}
\pf For all $e'\in\Gamma(E)$, by condition (ii) in Definition \ref{MTLA}, we have
\begin{eqnarray*}
\rho(e)\rho(e')(f)=\rho(e)(Df,e')_-&=&(e\star Df-\frac{1}{2}D(e,Df)_-,e')_-+(Df,[e,e']_E)_-\\
&=&(e\star Df,e')_-+\frac{1}{2}\rho(e')\rho(e)(f)+[\rho(e),\rho(e')](f),
\end{eqnarray*}
which implies that
$$(e\star Df,e')_-=(\frac{1}{2}D(Df,e)_-,e')_-.$$
Since the skew-symmetric bilinear form  $(\cdot,\cdot)_-$ is nondegenerate, we deduce that \eqref{pro1 MLA} holds.\qed\vspace{3mm}

 Given a symplectic Lie algebroid $(E,[\cdot,\cdot]_E,\rho,\omega)$, define a multiplication $\star:\Gamma(E)\times\Gamma(E)\longrightarrow \Gamma(E)$  by
\begin{eqnarray}\label{LSCA bracket}
\omega(e_1{\star}e_2,e_3)=\rho(e_1)\omega(e_2,e_3)+\frac{1}{2}\rho(e_3)\omega(e_1,e_2)-\omega(e_2,[e_1,e_3]_E),\quad \forall e_1,e_2,e_3\in \Gamma(E).
\end{eqnarray}
Or equivalently,
\begin{equation}\label{eq:newmul}
  e_1\star e_2={\omega^\sharp}^{-1}(\huaL_{e_1}\omega^\sharp(e_2)+\half \dM(\omega(e_1,e_2))).
\end{equation}
Actually, by \eqref{LSCA bracket}, we have
\begin{eqnarray*}
  \langle\omega^\sharp(e_1\star e_2),e_3\rangle&=&\rho(e_1)\langle\omega^\sharp(e_2),e_3\rangle+\half \langle \dM(\omega(e_1,e_2)),e_3\rangle-\langle\omega^\sharp(e_2),[e_1,e_3]_E\rangle\\
  &=&\langle \huaL_{e_1}\omega^\sharp(e_2),e_3\rangle+\half \langle \dM(\omega(e_1,e_2)),e_3\rangle.
\end{eqnarray*}
Therefore, we have
\begin{equation*}
  \omega^\sharp(e_1\star e_2)= \huaL_{e_1}\omega^\sharp(e_2)+\half \dM(\omega(e_1,e_2)),
\end{equation*}
which implies that \eqref{LSCA bracket} and \eqref{eq:newmul} are equivalent.

\begin{thm}\label{thm:Symplectic-LWX}
Let $(E,[\cdot,\cdot]_E,\rho,\omega)$ be a symplectic Lie algebroid. Then $(E,\star,\rho,(\cdot,\cdot)_-=\omega)$ is a pre-symplectic algebroid, and satisfies
\begin{eqnarray}\label{bracket}
[e_1,e_2]_E=e_1\star e_2-e_2\star e_1,\quad \forall e_1,e_2\in \Gamma(E),
\end{eqnarray}
where the multiplication $\star$ is given by $(\ref{LSCA bracket})$.
\end{thm}
\pf Since $\dM\omega(e_1,e_2,e_3)=0$, we have
\begin{eqnarray*}
&&\omega([e_1,e_2]_E,e_3)\\
&=&\rho(e_1)\omega(e_2,e_3)-\rho(e_2)\omega(e_1,e_3)+\rho(e_3)\omega(e_1,e_2)-\omega([e_2,e_3]_E,e_1)+\omega([e_1,e_3]_E,e_2).
\end{eqnarray*}
On the other hand, by (\ref{LSCA bracket}), we have
\begin{eqnarray*}
&&\omega(e_1\star e_2-e_2\star e_1,e_3)\\
&=&\rho(e_1)\omega(e_2,e_3)-\rho(e_2)\omega(e_1,e_3)+\rho(e_3)\omega(e_1,e_2)-\omega([e_2,e_3]_E,e_1)+\omega([e_1,e_3]_E,e_2).
\end{eqnarray*}
Since $\omega$ is nondegenerate, we have
$$[e_1,e_2]_E=e_1\star e_2-e_2\star e_1.$$

By (\ref{LSCA bracket}), it is straightforward to deduce that condition (ii) in Definition \ref{MTLA} holds.
 For all $e_1,e_2,e_3,e_4\in \Gamma(E)$, we have
\begin{eqnarray*}
([e_1,e_2]_E\star e_3,e_4)_-&=&\rho([e_1,e_2]_E)(e_3,e_4)_-+\frac{1}{2}\rho(e_4)([e_1,e_2]_E,e_3)_--(e_3,[[e_1,e_2]_E,e_4]_E)_-;\\
(e_1\star(e_2\star e_3)_,e_4)_-&=&\rho(e_1)\rho(e_2)(e_3,e_4)_-+\frac{1}{2}\rho(e_1)\rho(e_4)(e_2,e_3)_--\rho(e_1)(e_3,[e_2,e_4]_E)_-\\
&&-\frac{1}{2}\rho(e_4)\rho(e_2)(e_3,e_1)_-
-\frac{1}{4}\rho(e_4)\rho(e_1)(e_2,e_3)_-+\frac{1}{2}\rho(e_4)([e_1,e_2]_E,e_3)_-\\
&&+\rho(e_2)([e_1,e_4]_E,e_3)_--\frac{1}{2}\rho([e_1,e_4]_E)(e_2,e_3)_-+(e_3,[e_2,[e_1,e_4]_E]_E)_-;\\
(e_2\star(e_1\star e_3)_,e_4)_-&=&\rho(e_2)\rho(e_1)(e_3,e_4)_-+\frac{1}{2}\rho(e_2)\rho(e_4)(e_1,e_3)_--\rho(e_2)(e_3,[e_1,e_4]_E)_-\\
&&-\frac{1}{2}\rho(e_4)\rho(e_1)(e_3,e_2)_-
-\frac{1}{4}\rho(e_4)\rho(e_2)(e_1,e_3)_-+\frac{1}{2}\rho(e_4)([e_2,e_1]_E,e_3)_-\\
&&+\rho(e_1)([e_2,e_4]_E,e_3)_--\frac{1}{2}\rho([e_2,e_4]_E)(e_1,e_3)_-+(e_3,[e_1,[e_2,e_4]_E]_E)_-.
\end{eqnarray*}
Therefore, we obtain
\begin{eqnarray*}
&&\big((e_1,e_2,e_3)-(e_2,e_1,e_3),e_4\big)_-\\
&=&-([e_1,e_2]_E\star e_3,e_4)_-+(e_1\star(e_2\star e_3)_,e_4)_--(e_2\star(e_1\star e_3)_,e_4)_-\\
&=&\frac{1}{2}\rho(e_4)([e_1,e_2]_E,e_3)_-+\frac{1}{4}\rho(e_4)\rho(e_2)(e_1,e_3)_--\frac{1}{4}\rho(e_4)\rho(e_1)(e_2,e_3)_-.
\end{eqnarray*}
On the other hand, by \eqref{T-equation} and \eqref{LSCA bracket}, we have
\begin{eqnarray*}
T(e_1,e_2,e_3)&=&-\frac{1}{2}\rho(e_1)(e_2,e_3)_-+\frac{1}{2}\rho(e_2)(e_1,e_3)_-+\rho(e_3)(e_1,e_2)_-\\
&&+2([e_1,e_2]_E,e_3)_-+([e_1,e_3]_E,e_2)_--([e_2,e_3]_E,e_1)_-.
\end{eqnarray*}
By the fact that $(\cdot,\cdot)_-=\omega$ is closed, we have
\begin{eqnarray*}
&&-([e_2,e_3]_E,e_1)_-+([e_1,e_3]_E,e_2)_-\\
&=&([e_1,e_2],e_3)_--\rho(e_1)(e_2,e_3)_-+\rho(e_2)(e_1,e_3)_--\rho(e_3)(e_1,e_2)_-.
\end{eqnarray*}
Therefore, we have
\begin{eqnarray*}
T(e_1,e_2,e_3)=3([e_1,e_2]_E,e_3)_-+\frac{3}{2}\rho(e_2)(e_1,e_3)_--\frac{3}{2}\rho(e_1)(e_2,e_3)_-,
\end{eqnarray*}
which implies that
\begin{eqnarray*}
\frac{1}{6}(DT(e_1,e_2,e_3),e_4)_-&=&\frac{1}{2}\rho(e_4)([e_1,e_2]_E,e_3)_-+\frac{1}{4}\rho(e_4)\rho(e_2)(e_1,e_3)_--\frac{1}{4}\rho(e_4)\rho(e_1)(e_2,e_3)_-\\
&=&\big((e_1,e_2,e_3)-(e_2,e_1,e_3),e_4\big)_-.
\end{eqnarray*}
Thus, condition (i) in Definition \ref{MTLA} holds.  \qed\vspace{3mm}

By Theorem \ref{thm:Symplectic-LWX}, we can obtain many examples of pre-symplectic algebroids.

\begin{ex}{\rm
Let $(\huaA,[\cdot,\cdot]_\huaA,a_\huaA)$ be a Lie algebroid  and  $\Lambda\in\Gamma(\wedge^2\huaA)$ an invertible bisection of $\huaA$ satisfying $[\Lambda,\Lambda]=0$. Then $(\huaA^*,[\cdot,\cdot]_{\Lambda},a_{\huaA^*}=a_\huaA\circ\Lambda^\sharp,\Lambda)$ is a symplectic Lie algebroid, where the bundle map $\Lambda^\sharp:\huaA^*\longrightarrow \huaA $ is defined by
\begin{equation*}
\Lambda^\sharp(\xi)(\eta)=\Lambda(\xi,\eta),\quad \forall \xi,\eta\in\Gamma(\huaA^*),
\end{equation*}
and the bracket $[\cdot,\cdot]_{\Lambda}:\wedge^2\Gamma(\huaA^*)\longrightarrow\Gamma(\huaA^*)$ is defined by
\begin{equation*}
[\xi,\eta]_{\Lambda}=\huaL_{\Lambda^\sharp(\xi)}\eta-\huaL_{\Lambda^\sharp(\eta)}\xi-\dM \Lambda(\xi,\eta),\quad \forall\xi,\eta\in\Gamma(\huaA^*).
\end{equation*}
 By \eqref{eq:newmul}, and the fact $$(\Lambda^\sharp)^{-1}(\huaL_\xi\Lambda^\sharp(\eta))=\huaL_{\Lambda^\sharp(\xi)}\eta,\quad (\Lambda^\sharp)^{-1}(\dM_*f)=-\dM f,$$ the pre-symplectic algebroid structure on $\huaA^*$ is given by
\begin{equation}
\xi{\star}\eta=\huaL_{\Lambda^\sharp(\xi)}\eta-\half\dM\Lambda(\xi,\eta),\quad \forall \xi,\eta\in \Gamma(\huaA^*).
\end{equation}

}
\end{ex}

\begin{ex}{\rm
Consider  $\mathbb R^{2n}$ and the standard symplectic structure $\omega=\sum_{i=1}^nd x_i\wedge dx_{n+i}$. Then $(T\mathbb R^{2n},\omega)$ is a symplectic Lie algebroid. By \eqref{eq:newmul}, the pre-symplectic algebroid structure is given by
\begin{eqnarray*}
  \frac{\partial}{\partial x_l}\star  \frac{\partial}{\partial x_{m}}&=&0,\quad 1\leq l,m\leq 2n,\\
   \frac{\partial}{\partial x_l}\star  f\frac{\partial}{\partial x_{m}}&=&\frac{\partial f}{\partial x_{l}}\frac{\partial}{\partial x_m},~~(f\frac{\partial}{\partial x_l})\star  \frac{\partial}{\partial x_{m}}=0,\quad 1\leq l,m\leq 2n,\mid l-m\mid\neq n,\\
    \frac{\partial}{\partial x_i}\star  f\frac{\partial}{\partial x_{n+i}}&=&\half\Big(\frac{\partial f}{\partial x_i}\frac{\partial}{\partial x_{n+i}}+\frac{\partial f}{\partial x_{n+i}}\frac{\partial}{\partial x_i}\Big),\\
    (f\frac{\partial}{\partial x_{n+i}})\star \frac{\partial}{\partial x_i}&=&-\half\Big(\frac{\partial f}{\partial x_i}\frac{\partial}{\partial x_{n+i}}-\frac{\partial f}{\partial x_{n+i}}\frac{\partial}{\partial x_i}\Big),\\
     (f \frac{\partial}{\partial x_i})\star  \frac{\partial}{\partial x_{n+i}}&=&\half\Big(\frac{\partial f}{\partial x_i}\frac{\partial}{\partial x_{n+i}}-\frac{\partial f}{\partial x_{n+i}}\frac{\partial}{\partial x_i}\Big),\\
     \frac{\partial}{\partial x_{n+i}}\star f\frac{\partial}{\partial x_i}&=&\half\Big(\frac{\partial f}{\partial x_i}\frac{\partial}{\partial x_{n+i}}+\frac{\partial f}{\partial x_{n+i}}\frac{\partial}{\partial x_i}\Big).
\end{eqnarray*}
}
\end{ex}

The following examples is due to the fact  that
for any regular Poisson manifold $(M,\pi)$, its characteristic
distribution $\Img(\pi^\sharp)\subset TM$ is a symplectic Lie
algebroid, where $\pi^\sharp:T^*M\longrightarrow TM$ is given by
$\langle\pi^\sharp(\xi),\eta\rangle=\pi(\xi,\eta)$ for all
$\xi,\eta\in\Omega^1(M)$, and the symplectic form is characterized
by $\omega_p(H_f,H_g)=\pi(df,dg)|_p$, in which $p\in M,~f, g\in
C^\infty(M)$ and $H_f, H_g$ are Hamiltonian vector fields.
\begin{ex}{\rm
 Consider the Poisson manifold $(\mathbb R^3/\{0\},\pi=x\frac{\partial}{\partial y}\wedge \frac{\partial}{\partial z}+y\frac{\partial}{\partial z}\wedge \frac{\partial}{\partial x}+z\frac{\partial}{\partial x}\wedge \frac{\partial}{\partial y}).$  Then $\Img(\pi^\sharp)$ is a symplectic Lie algebroid. In the sequel, we give the corresponding pre-symplectic algebroid structure locally. At any point $p=(x,y,z)$, where $y\neq 0$, $$\{e_1=y \frac{\partial}{\partial x}-x\frac{\partial}{\partial y},~e_2=z \frac{\partial}{\partial y}-y\frac{\partial}{\partial z}\}$$ is a basis of $\Img(\pi_p^\sharp)$.
We know  that symplectic leaves are of the form
$$\{(x,y,z)\in\mathbb R^3:x^2+y^2+z^2=r^2,r\neq0\}.$$
On every symplectic leaf, the symplectic structure is given by
$
\omega=ye_1^*\wedge e_2^*.
$
By \eqref{eq:newmul}, we have
\begin{eqnarray*}
e_1\star e_2=-\frac{1}{2y}(ze_1+xe_2);\quad
e_2\star e_1=\frac{1}{2y}(ze_1+xe_2).
\end{eqnarray*}
}
\end{ex}

The following example was introduced in \cite{symplectic lie algebroid} which plays an important role in geometric mechanics.
\begin{ex}{\rm
Let $(A,[\cdot,\cdot]_A,a)$ be a Lie algebroid over a manifold $M$, and $\tau:A\longrightarrow M$   the vector bundle projection. Let   $\tau^*:A^*\longrightarrow M$ be the vector bundle projection of the dual bundle $A^*$. We consider the prolongation $\huaL^{\tau^*}A$ of $A$ over $\tau^*$,
 $$
 \huaL^{\tau^*}A=\{(x,v)\in A\times TA^*\mid a(x)=T\tau^*(v)\},
 $$
 where $T\tau^*:TA^*\longrightarrow TM$ is the tangent map of $\tau^*$. Denote by $\tau^{\tau^*}:\huaL^{\tau^*}A\longrightarrow A^*$ the map given by
$$
\tau^{\tau^*}(x,v):=\tau_{A^*}(v),\quad \forall (x,v)\in \huaL^{\tau^*}A,
$$
where $\tau_{A^*}:TA^*\longrightarrow A^*$ is the canonical projection. Now let $\{x^i\}$ be local coordinates on an open subset $U$ of $M$, $\{e_\alpha\}$
 a basis of sections of the vector bundle $\tau^{-1}(U)\longrightarrow U$ and $\{e^\alpha\}$  the dual basis of $\{e_\alpha\}$. Then $\{\tilde{e}_\alpha,\bar{e}_\alpha\}$ is
 the basis of the vector bundle $\tau^{\tau^*}((\tau^*)^{-1}(U))\longrightarrow(\tau^*)^{-1}(U)$, which is given by
 \begin{eqnarray*}
 \tilde{e}_\alpha(\xi)&=&(e_\alpha(\tau^*(\xi)),a_\alpha^i\frac{\partial}{\partial x^i}\mid_\xi),\\
 \bar{e}_\alpha(\xi)&=&(0,\frac{\partial}{\partial y_\alpha}\mid_\xi).
 \end{eqnarray*}
 Here $a_\alpha^i$ are the components of the anchor map with respect to the basis $\{e_\alpha\}$, i.e. $a(e_\alpha)={a}_\alpha^i\frac{\partial}{\partial x^i}$ and $\{x^i,y_\alpha\}$ are the corresponding local coordinates on $A^*$ induced by the local coordinates $\{x^i\}$ and the basis $\{e_\alpha\}$. Then the Lie algebroid structure $(\huaL^{\tau^*}A,[\cdot,\cdot]^{\tau^*},a^{\tau^*})$ on $\huaL^{\tau^*}A$ is locally defined by
 \begin{eqnarray*}
  [\tilde{e}_\alpha,\tilde{e}_\beta]^{\tau^*}&=&C_{\alpha \beta}^\gamma \tilde{e}_\gamma,\\
   {[\tilde{e}_\alpha,\bar{e}_\beta]}^{\tau^*}&=& [\bar{e}_\alpha,\bar{e}_\beta]^{\tau^*}=0,\\
  a^{\tau^*}(\tilde{e}_\alpha)&=&a_\alpha^i\frac{\partial}{\partial x^i},\quad  a^{\tau^*}(\bar{e}_\alpha)=\frac{\partial}{\partial y_\alpha},
   \end{eqnarray*}
    where $C_{\alpha\beta}^\gamma$ are structure functions of the Lie bracket $[\cdot,\cdot]_A$ with respect to the basis $\{e_\alpha\}$.

The canonical symplectic section $\Omega$ is locally given by
 \begin{equation}
 \Omega=\tilde{e}^\alpha\wedge\bar{e}^\alpha+\frac{1}{2}C_{\alpha\beta}^\gamma y_\gamma\tilde{e}^\alpha\wedge\tilde{e}^\beta.
 \end{equation}
 By \eqref{eq:newmul}, the pre-symplectic algebroid structure is given by
 \begin{eqnarray*}
 \tilde{e}_k\star\tilde{e}_k&=&(-{a}_k^iy_\gamma \frac{\partial C_{km}^\gamma}{\partial x^i}+C_{km}^\gamma C_{kr}^py_p)\bar{e}_m;\\
 \tilde{e}_k\star\tilde{e}_l&=&(C_{km}^\gamma C_{lr}^py_p-{a}_k^iy_\gamma \frac{\partial C_{lm}^\gamma}{\partial x^i}-\half{a}_m^iy_\gamma \frac{\partial C_{kl}^\gamma}{\partial x^i} +\half C_{kl}^p C_{jm}^\gamma y_\gamma)\bar{e}_m+\half C_{kl}^m\tilde{e}_m,\quad k\neq l;\\
 \tilde{e}_k\star\bar{e}_l&=&-C_{km}^l\bar{e}_m;\quad \quad\bar{e}_l\star\tilde{e}_k=-C_{km}^l\bar{e}_m;\quad \quad \bar{e}_k\star\bar{e}_l=0.
 \end{eqnarray*}
 }
\end{ex}
\vspace{3mm}

  In \cite{LiuShengBaiChen}, we proved that there is a symplectic Lie algebroid $A^c\ltimes_{L^*} A^*$ associated to any left-symmetric algebroid $A$ (Theorem \ref{thm:symLie}). Now it is natural to investigate the corresponding pre-symplectic algebroid.

\begin{pro}\label{pro:standard}
  Let $(A,\cdot_A,a_A)$ be a left-symmetric algebroid and $(A^c\ltimes_{L^*} A^*,\omega)$ the corresponding symplectic Lie algebroid, where $\omega$ is given by \eqref{eq:defiomega}. Then the corresponding pre-symplectic algebroid structure is given by
  \begin{equation}\label{eq:standardLWXbracket}
    (x+\xi)\star(y+\eta)=x\cdot_A y+\huaL_x\eta-R_y^*\xi-\half \dM(x+\xi,y+\eta)_+,
  \end{equation}
where $\huaL$ and $R^*$ are given by \eqref{Lie der1} and \eqref{eq:dualR} respectively and $(\cdot,\cdot)_+$ is the nondegenerate symmetric bilinear form on $A\oplus A^*$ given by
\begin{equation}\label{sym-form}
 (x+\xi,y+\eta)_+=\langle\xi,y\rangle+\langle\eta,x\rangle.
\end{equation}
\end{pro}
  We call the above pre-symplectic algebroid the {\bf pseudo-semidirect product} of the left-symmetric algebroid  $(A,\cdot_A,a_A)$ and its dual bundle $A^*$.

\pf By \eqref{LSCA bracket}, for all $x,y,z\in\Gamma(A),\xi,\eta,\zeta\in\Gamma(A^*)$, we have
\begin{eqnarray*}
&&\omega((x+\xi)\star (y+\eta),z+\zeta)\\
&=&\rho(x+\xi)\omega(y+\eta,z+\zeta)+\frac{1}{2}\rho(z+\zeta)\omega(x+\xi,y+\eta)-\omega(y+\eta,[x+\xi,z+\zeta]_S)\\
&=&a_A(x)(\langle \eta,z\rangle-\langle y,\zeta\rangle)+\frac{1}{2}a_A(z)(\langle \xi,y\rangle-\langle x,\eta\rangle)
-\langle \eta,[x,z]_A\rangle+\langle y,L^*_x\zeta-L^*_z\xi\rangle\\
&=&a_A(x)\langle \eta,z\rangle-\langle \eta,[x,z]_A\rangle-a_A(x)\langle y,\zeta\rangle+\langle y,L^*_x\zeta\rangle+\frac{1}{2}a_A(z)\langle \xi,y\rangle-\langle y,L^*_z\xi\rangle-\frac{1}{2}a_A(z)\langle x,\eta\rangle\\
&=&\langle\huaL_x\eta,z\rangle-\langle x\cdot_A y,\zeta\rangle-\langle R_y^*\xi,z\rangle-\frac{1}{2}a_A(z)(\langle\xi,y\rangle+\langle x,\eta\rangle)\\
&=&\omega(x\cdot_A y+\huaL_x\eta-R_y^*\xi-\half \dM(x+\xi,y+\eta)_+,z+\zeta),
\end{eqnarray*}
which implies that \eqref{eq:standardLWXbracket} holds.\qed

 \begin{defi}
 Let $(E,\star,\rho,(\cdot,\cdot)_-)$ be a pre-symplectic algebroid. A subbundle $F$ of $E$ is called {\bf isotropic}
  if it is isotropic under the skew-symmetric bilinear form $(\cdot,\cdot)_-$. It is called {\bf integrable}
  if $\Gamma(F)$ is closed under the operation $\star$. A {\bf Dirac structure} is a subbundle $F$ which is maximally isotropic and integrable.
 \end{defi}

The following proposition is obvious.
\begin{pro}\label{Dirac subbundles}
Let $F$ be a Dirac structure of a pre-symplectic algebroid $(E,\star,\rho,(\cdot,\cdot)_-)$. Then $(F,\star|_F,\rho|_F)$ is a left-symmetric algebroid.
\end{pro}

\begin{pro}\label{pro:Dirac-Lag}
There is a one-to-one correspondence between Dirac structures of a pre-symplectic algebroid $(E,\star,\rho,(\cdot,\cdot)_-)$ and Lagrangian subalgebroids of  the corresponding symplectic Lie algebroid $(E,[\cdot,\cdot]_E,\rho,\omega=(\cdot,\cdot)_-)$.
\end{pro}
\pf Let $F$ be a Dirac structure of $(E,\star,\rho,(\cdot,\cdot)_-)$. By Proposition \ref{Dirac subbundles}, $(F,\star|_F,\rho|_F)$ is a left-symmetric algebroid.  Then it is straightforward to see that $(F,[\cdot,\cdot]_F={[\cdot,\cdot]_E}|_F,\rho|_F)$ is a Lagrangian subalgebroid of  the corresponding symplectic Lie algebroid.

Conversely, assume that $(F,[\cdot,\cdot]_E|_F,\rho|_F)$ is a Lagrangian subalgebroid of  the corresponding symplectic Lie algebroid. By \eqref{bracket}, for all $e_1,e_2,e_3\in\Gamma(F)=\Gamma(F^\bot)$, we have
\begin{eqnarray*}
\omega(e_1{\star}e_2,e_3)=\rho(e_1)\omega(e_2,e_3)+\frac{1}{2}\rho(e_3)\omega(e_1,e_2)-\omega(e_2,[e_1,e_3]_E)=0,
\end{eqnarray*}
which implies that $e_1{\star}e_2\in\Gamma(F)$, i.e. $\Gamma(F)$ is closed under the operation $\star$. Therefore, $F$ is a Dirac structure.\qed

\section{Classification of exact  pre-symplectic algebroids}

In this section, we study  exact  pre-symplectic algebroids, which
can be viewed as twists of the pseudo-semidirect products by
3-cocycles. In the sequel, $(M,\nabla)$ is a flat manifold and $T_\nabla M$ is the associated
left-symmetric algebroid given in Example \ref{ex:main}.

\begin{defi}
Let $(M,\nabla)$ be a flat manifold. A
pre-symplectic algebroid $(E,\star,\rho,(\cdot,\cdot)_-)$ is called
{\bf exact with respect to the left-symmetric algebroid $T_\nabla M$}
if the following short sequence of vector bundles
$$ 0\longrightarrow T^*M\stackrel{\rho^*}{\longrightarrow}E\stackrel{\rho}\longrightarrow TM\longrightarrow0$$
is exact and $\rho(e_1\star e_2)=\nabla_{\rho(e_1)}\rho(e_2)$, where  $\rho^*:T^*M\longrightarrow E$ is given by
\begin{eqnarray*}
( \rho^*(\xi),e)_-=\langle\xi,\rho(e)\rangle,\quad \forall \xi\in \Omega^1(M),~~e\in \Gamma(E).
\end{eqnarray*}
\end{defi}
Now choose an isotropic splitting $\sigma:TM\longrightarrow E$ of the above exact sequence.  Define $\phi:\XM \times \XM\longrightarrow \Omega^1(M)$ by
\begin{equation}\label{eq:3-form}
\rho^*\phi(x,y)=\sigma(x)\star\sigma(y)-\sigma(\nabla_xy),\quad \forall x,y\in\XM.
\end{equation}
By Proposition \ref{Condii iii}, we have $\phi\in \Gamma(\otimes^3T^*M)$. In fact, by \eqref{eq:3-form}, we obtain
\begin{eqnarray*}
\rho^*\phi(x,fy)&=&\sigma(x)\star\sigma(fy)-\sigma(\nabla_x fy)\\
&=&f\sigma(x)\star\sigma(y)+\rho\circ\sigma(x)(f)y-f\sigma(\nabla_xy)-x(f)\sigma(y)\\
&=&f\rho^*\phi(x,y).
\end{eqnarray*}
Similarly, we have
$$\rho^*\phi(fx,y)=f\rho^*\phi(x,y).$$
Given an isotropic splitting $\sigma$, we have $E\cong TM\oplus T^*M$. Under this isomorphism, the nondegenerate skew-symmetric bilinear form $(\cdot,\cdot)_- $
is exactly given by \eqref{eq:defiomega},  $\rho$ is the projection to the first summand. Using condition (ii) in Definition \ref{MTLA}, we can deduce that the multiplication $\star$ is given by\footnote{It is a twist of the pseudo-semidirect product given by \eqref{eq:standardLWXbracket}. }
 \begin{equation}\label{eq:3-formbracket}
    (x+\xi)\star(y+\eta)=\nabla_xy+\huaL_x\eta-R_y^*\xi-\half \dM(x+\xi,y+\eta)_++\phi(x,y).
  \end{equation}

By condition (ii) in Definition \ref{MTLA}, we obtain
\begin{equation}\label{eq:3-form ii}
\langle\phi(x,y),z\rangle=\langle y,\phi(x,z)-\phi(z,x)\rangle,
\end{equation}
which implies that
 \begin{equation}\label{eq:3-form ii2}
\phi(x,y,z)=\phi(x,z,y)-\phi(z,x,y).
\end{equation}
Therefore, we have
\begin{eqnarray}\label{eq:13skewsym}
\phi(x,y,z)=-\phi(z,y,x).
\end{eqnarray}

Now we define $\tilde{\phi}\in\Gamma(\wedge^2T^*M\otimes T^*M)$ by
 \begin{equation}\label{eq:3-form2}
 \tilde{\phi}(x,y,z)=\phi(x,z,y),\quad \forall x,y,z\in\XM.
 \end{equation}
By \eqref{eq:3-form ii2} and \eqref{eq:13skewsym}, we have
 \begin{equation}
 \tilde{\phi}(x,y,z)+c.p.=0,\quad \forall x,y,z\in\XM.
 \end{equation}

\begin{pro}With the above notations, $(TM\oplus T^*M,\star,\rho,(\cdot,\cdot)_-)$ is a pre-symplectic algebroid, where the nondegenerate skew-symmetric bilinear form $(\cdot,\cdot)_- $
is   given by \eqref{eq:defiomega},  $\rho$ is the projection to the first summand and the multiplication $\star$ is given by \eqref{eq:3-formbracket}, if and only if $\tilde{\phi} $ is closed, i.e.
\begin{eqnarray}
\delta\tilde{\phi}=0.
\end{eqnarray}
Here $\delta$ is the coboundary operator associated to the left-symmetric algebroid $T_\nabla M$ given by \eqref{LSCA cohomology}.
\end{pro}
\pf First it is straightforward to deduce that condition (ii) in Definition \ref{MTLA} holds. Then by \eqref{eq:3-form ii2}, for all $x,y,z\in\XM$ and $\xi,\eta,\zeta\in\Omega^1(M)$, we have
\begin{eqnarray*}
&&(x+\xi,y+\eta,z+\zeta)-(y+\eta,x+\xi,z+\zeta)-\frac{1}{6}DT(x+\xi,y+\eta,z+\zeta)\\
&=&\phi(x,\nabla_yz)-\phi(\nabla_xy,z)+\huaL_{x}\phi(y,z)+R_z^*\phi(x,y)\\
&&-\phi(y,\nabla_xz)+\phi(\nabla_yx,z)-\huaL_{y}\phi(x,z)-R_z^*\phi(y,x)\\
&&+\frac{1}{3}\dM(\phi(x,y,z)-\phi(y,z,x)-\phi(y,x,z)+\phi(x,z,y))\\
&=&\phi(x,\nabla_yz)-\phi(\nabla_xy,z)+\huaL_{x}\phi(y,z)+R_z^*\phi(x,y)\\
&&-\phi(y,\nabla_xz)+\phi(\nabla_yx,z)-\huaL_{y}\phi(x,z)-R_z^*\phi(y,x)+\dM(\phi(x,z,y)),
\end{eqnarray*}
and
\begin{eqnarray*}
&&\langle\phi(x,\nabla_yz)-\phi(\nabla_xy,z)+\huaL_{x}\phi(y,z)+R_z^*\phi(x,y)\\
&&-\phi(y,\nabla_xz)+\phi(\nabla_yx,z)-\huaL_{y}\phi(x,z)-R_z^*\phi(y,x)+\dM(\phi(x,z,y)),w\rangle\\
&=&\phi(x,\nabla_yz,w)-\phi(\nabla_xy,z,w)+x\phi(y,z,w)-\phi(y,z,[x,w])\\
&&-\phi(x,y,\nabla_wz)-\phi(y,\nabla_xz,w)+\phi(\nabla_yx,z,w)-y\phi(x,z,w)\\
&&+\phi(x,z,[y,w])+\phi(y,x,\nabla_wz)+w\phi(x,z,y)\\
&=&\phi(x,\nabla_yz,w)-\phi([x,y],z,w)+x\phi(y,z,w)-\phi(y,z,[x,w])\\
&&-\phi(x,\nabla_wz,y)-\phi(y,\nabla_xz,w)-y\phi(x,z,w)+\phi(x,z,[y,w])+w\phi(x,z,y)\\
&=&\tilde{\phi}(x,w,\nabla_yz)-\tilde{\phi}([x,y],w,z)+x\tilde{\phi}(y,w,z)-\tilde{\phi}(y,[x,w],z)\\
&&-\tilde{\phi}(x,y,\nabla_wz)-\tilde{\phi}(y,w,\nabla_xz)-y\tilde{\phi}(x,w,z)+\tilde{\phi}(x,[y,w],z)+w\tilde{\phi}(x,y,z)\\
&=&\delta\tilde{\phi}(x,y,w,z).
\end{eqnarray*}
Therefore, condition (i) in Definition \ref{MTLA} holds if and only if $\tilde{\phi}$ is closed. \qed\vspace{3mm}

Let $\sigma':TM\longrightarrow E$ be another isotropic splitting.  Assume that
\begin{equation}\label{eq:3-coboundary}
\sigma(x)-\sigma'(x)=\theta(x),\quad \forall x\in\XM,
\end{equation}
where $\theta:TM\longrightarrow T^*M$ is a bundle map. Similar to \eqref{eq:3-form}, define $\phi':\XM \times \XM\longrightarrow \Omega^1(M)$ by
$$
\rho^*\phi'(x,y)=\sigma'(x)\star\sigma'(y)-\sigma'(\nabla_xy),\quad \forall x,y\in\XM.
$$
Since both $\sigma$ and $\sigma'$ are isotropic, we have
\begin{eqnarray*}
0&=&(\sigma(x),\sigma(y))_-=(\sigma'(x)+\rho^*\theta(x),\sigma'(y)+\rho^*\theta(y))_-\\
&=&(\rho^*\theta(x),\sigma'(y))_--(\rho^*\theta(y),\sigma'(x))_-\\
&=&\theta(x,y)-\theta(y,x),
\end{eqnarray*}
which implies that $$\theta(x,y)=\theta(y,x).$$

Furthermore, by \eqref{eq:3-form} and condition (ii) in Definition \ref{MTLA}, we have
\begin{eqnarray*}
\rho^*(\phi(x,y)-\phi'(x,y))&=&\sigma(x)\star\sigma(y)-\sigma'(x)\star\sigma'(y)-\big(\sigma(\nabla_xy)-\sigma'(\nabla_xy)\big)\\
&=&\sigma'(x)\star \rho^*\theta(y)+\rho^*\theta(x)\star\sigma'(y)-\rho^*\theta(\nabla_xy)\\
&=&\rho^*\Big(\huaL_x\theta(y)-R_y^*\theta(x)-\frac{1}{2}\dM(\theta(y,x)+\theta(x,y))-\theta(\nabla_xy)\Big).
\end{eqnarray*}
Since $\rho$ is surjective, $\rho^*$ is injective. Therefore, we have
\begin{eqnarray*}
\phi(x,y,z)-\phi'(x,y,z)&=&\langle\huaL_x\theta(y)-R_y^*\theta(x)-\frac{1}{2}\dM(\theta(y,x)+\theta(x,y))-\theta(\nabla_xy),z\rangle\\
&=&x\theta(y,z)-\theta(y,[x,z])-z\theta(x,y)+\theta(x,\nabla_zy)-\theta(\nabla_xy,z)\\
&=&\delta\theta(x,z,y).
\end{eqnarray*}
Equivalently, we have
\begin{equation}\label{eq:diffexact}
\tilde{\phi}(x,z,y)-\tilde{\phi'}(x,z,y)=\delta\theta(x,z,y).
\end{equation}
Now we introduce a new cochain complex, whose third cohomology group classifies exact pre-symplectic algebroids. Let $(A,\cdot_A,a_A)$ be a left-symmetric algebroid. For all $x,y,z\in \Gamma(A)$, define
\begin{eqnarray*}
\tilde{C}^1(A)&=&\{\varphi\in C^1(A)\mid \varphi([x,y]_A)=a_A(x)\varphi(y)-a_A(y)\varphi(x)\},\\
\tilde{C}^2(A)&=&\{\varphi\in C^2(A)\mid \varphi(x,y)=\varphi(y,x)\},\\
\tilde{C}^3(A)&=&\{\varphi\in C^3(A)\mid \varphi(x,y,z)+c.p.=0\},\\
\tilde{C}^n(A)&=&C^n(A),\quad n\geq 4.
 \end{eqnarray*}
It is straightforward to verify that the cochain complex $(\tilde{C}^{n}(A),\delta)$ is a subcomplex of the cochain complex $(C^{n}(A),\delta)$. We denote the $n$-th cohomology group of $(\tilde{C}^{n}(A),\delta)$ by $\tilde{H}^n(A)$.

By \eqref{eq:diffexact}, we have
\begin{thm} Let  $(M,\nabla)$ be a flat manifold and $T_\nabla M$   the associated
left-symmetric algebroid.
Then
the  exact pre-symplectic algebroids with respect to the left-symmetric algebroid $T_\nabla M$ are classified by the third  cohomology group $\tilde{H}^3(T_\nabla M )$.
\end{thm}

\section{Para-complex  pre-symplectic algebroids}

In this section, we give an important application of
pre-symplectic algebroids in the theory of pseudo-Riemannian Lie
algebroids. We introduce the notion of a  para-complex
pre-symplectic algebroid $(E,\star,\rho,(\cdot,\cdot)_-,P)$, which
can give rise to a para-K\"{a}hler Lie algebroid
$(E,[\cdot,\cdot]_E,\rho,\omega,P)$ and  a pseudo-Riemannian Lie
algebroid $(E,[\cdot,\cdot]_E,\rho,g)$. Thus, para-complex
pre-symplectic algebroids can be viewed as the underlying
structures of para-K\"{a}hler Lie algebroids. Moreover, the
multiplication $\star$ in a para-complex pre-symplectic algebroid
can characterize the Levi-Civita $E$-connection on Lagrangian
subalgebroids $E^+$ and $E^-$. This gives a geometric
interpretation of the multiplication  in a para-complex
pre-symplectic algebroid. Furthermore, there is a natural example
of para-complex pre-symplectic algebroid
constructed
 from
a pseudo-semidirect product.

\subsection{Para-K\"{a}hler and pseudo-Riemannian Lie algebroids}
\begin{defi}
 A {\bf
para-complex  structure} on  a real Lie algebroid $(E,[\cdot,\cdot]_E,\rho)$ over $M$ is a bundle map  $P:E\longrightarrow
E$ satisfying $P^2=\id$
 and the following integrability  condition:
\begin{equation}\label{eq:int-cond1}
P[e_1,e_2]_E=[P(e_1),e_2]_E+[e_1,P(e_2)]_E-P[P(e_1),P(e_2)]_E,\quad\forall e_1,e_2\in\Gamma(E).
\end{equation}
\end{defi}

We denote
\begin{eqnarray*}
E^+_m:=\{x\in E_m\mid P_m(x)=x,~\forall m\in M\},\quad E^-_m:=\{\xi\in E_m\mid P_m(\xi)=-\xi,~\forall m\in M\}.
\end{eqnarray*}
It is obvious that $E^+$ and $E^-$ are subbundles of  $E$  and $E$ is the Whitney sum of $E^+$ and $E^-$, i.e. $E=E^+\oplus E^-$.

\begin{pro}\label{Pro:int-equi}
The integrability  condition of $P$ is equivalent to
the fact
that $(E^+,[\cdot,\cdot{]\mid}_{E^+},{\rho\mid}_{E^+})$ and
$(E^-,[\cdot,\cdot]\mid_{E^-},\rho\mid_{E^-})$ are subalgebroids of
the Lie algebroid $(E,[\cdot,\cdot]_E,\rho)$.
\end{pro}

\begin{defi} \label{defi:parakahler}
A {\bf para-K\"{a}hler Lie algebroid} is a symplectic Lie algebroid $(E,[\cdot,\cdot]_E,\rho,\omega)$ with a para-complex structure $P$ satisfying the following condition:
\begin{equation}\label{eq:Para-Comp}
\omega(P(e_1),e_2)+\omega(e_1,P(e_2))=0,\quad \forall e_1,e_2\in \Gamma(E).
\end{equation}
\end{defi}
Denote a para-K\"{a}hler Lie algebroid by $(E,[\cdot,\cdot]_E,\rho,\omega,P)$ or simply by $(E,\omega,P)$.

 \begin{rmk}
The notion of a para-K\"{a}hler Lie algebroid was originally  introduced in \cite{Para-kahler}, in which a para-K\"{a}hler Lie algebroid is a symplectic Lie algebroid with two transversal Lagrangian subalgebroids. Our definition of a para-K\"{a}hler Lie algebroid is equivalent to the one given in \cite{Para-kahler}.

In fact, if $A_1$ and $A_2$ are two transversal Lagrangian subalgebroids of a symplectic Lie algebroid $(E,\omega)$, define $P:E\longrightarrow E$ by
\begin{equation}
P(x+\xi)=x-\xi,\quad x\in\Gamma(A_1),\xi\in\Gamma(A_2).
\end{equation}
Then $(E,\omega,P)$ is a para-K\"{a}hler Lie algebroid in the sense of Definition \ref{defi:parakahler}.

Conversely, let $(E,\omega,P)$ be a para-K\"{a}hler Lie algebroid in the sense of Definition \ref{defi:parakahler}. By \eqref{eq:Para-Comp}, both $(E^+,[\cdot,\cdot{]\mid}_{E^+},{\rho\mid}_{E^+})$ and $(E^-,[\cdot,\cdot]\mid_{E^-},\rho\mid_{E^-})$   are isotropic subalgebroids.  Since $\Gamma(E)=\Gamma(E^+)\oplus\Gamma(E^-)$,   $(E^+,[\cdot,\cdot{]\mid}_{E^+},{\rho\mid}_{E^+})$ and $(E^-,[\cdot,\cdot]\mid_{E^-},\rho\mid_{E^-})$ are Lagrangian subalgebroids of the symplectic Lie algebroid $(E,\omega)$.
 \end{rmk}

Recall that a {\bf pseudo-Riemannian Lie algebroid} is a Lie algebroid $(E,[\cdot,\cdot]_E,\rho)$ with a  pseudo-Riemannian metric $g$.

 \begin{defi}\label{defi:LCconnection}
 A connection $\nabla^g$  is called a  {\bf Levi-Civita $E$-connection} associated to the pseudo-Riemannian metric $g$ if for all $e_1,e_2,e_3\in\Gamma(E)$,  the following two conditions are satisfied:
\begin{itemize}
\item[$\rm(i)$] $\nabla^g$ is metric. i.e. $\rho(e_1)g(e_2,e_3)=g(\nabla^g_{e_1}e_2,e_3)+g(e_2,\nabla^g_{e_1}e_3);$
\item[$\rm(ii)$] $\nabla^g$ is torsion free. i.e. $[e_1,e_2]_E=\nabla^g_{e_1}e_2-\nabla^g_{e_2}e_1.$
\end{itemize}
\end{defi}
The Levi-Civita $E$-connection $\nabla^g$ associated to the pseudo-Riemannian metric $g$ can be uniquely determined by the following formula
\begin{eqnarray}
\nonumber2g(\nabla^g_{e_1}e_2,e_3)&=&\rho(e_1)g(e_2,e_3)+\rho(e_2)g(e_1,e_3)-\rho(e_3)g(e_1,e_2)\\
\label{eq:nablag}&&+g([e_3,e_1]_E,e_2)+g([e_3,e_2]_E,e_1)+g( [e_1,e_2]_E,e_3),~~ \forall e_1,e_2,e_3\in\Gamma(E).
\end{eqnarray}

Let $(E,\omega,P)$ be a para-K\"{a}hler Lie algebroid. Define $g\in\otimes^2 E^*$ by
\begin{equation}\label{eq:ParaC-Rie}
g(e_1,e_2):=\omega(e_1,P(e_2)),\quad \forall e_1,e_2\in\Gamma(E).
\end{equation}
By   \eqref{eq:Para-Comp}, $g$ is symmetric. By the non-degeneracy of $\omega$,  $g$ is non-degenerate. Therefore, $g$ is a pseudo-Riemannian metric on $A$.

\begin{pro}\label{pro:Levi-PC}
Let $(E,\omega,P)$ be a para-K\"{a}hler Lie algebroid. Then $(E,g)$ is a pseudo-Riemannian Lie algebroid, where $g$ is given by \eqref{eq:ParaC-Rie}. Furthermore, for all $e_1,e_2\in\Gamma(E)$, we have
\begin{eqnarray}
\nabla^g_{e_1}P(e_2)&=&P(\nabla^g_{e_1} e_2),\label{eq:Levi-PC}\\
g(P(e_1),P(e_2))&=&-g(e_1,e_2).\label{eq:Rie-PC}
\end{eqnarray}
\end{pro}
\pf By \eqref{eq:Para-Comp} and \eqref{eq:ParaC-Rie}, it is straightforward to deduce that \eqref{eq:Rie-PC} holds. Now we prove that \eqref{eq:Levi-PC} holds. For all $e\in\Gamma(E),y,z\in\Gamma(E^+)$, since $\Gamma(E^+)$ is a Lagrangian subalgebroid and $\omega$ is closed, we have
\begin{eqnarray*}
2g(\nabla^g_ey,z)&=&\rho(e)g(y,z)+\rho(y)g(e,z)-\rho(z)g(e,y)\\
&&+g([z,e]_E,y)+g([z,y]_E,e)+g( [e,y]_E,z)\\
&=&-\rho(e)\omega(y,z)+\rho(y)\omega(e,z)-\rho(z)\omega(e,y)\\
&&+\omega([z,e]_E,y)+\omega([y,z]_E,e)+\omega( [e,y]_E,z)\\
&=&-\dM \omega(e,y,z)=0,
\end{eqnarray*}
which implies that $g(\nabla^g_e y,z)=0$. Thus, for all $e\in\Gamma(E),y,z\in\Gamma(E^+)$, we have
\begin{eqnarray*}
g(\nabla^g_eP(y),z)&=&g(\nabla^g_ey,z)=0.
\end{eqnarray*}
On the other hand, by \eqref{eq:Para-Comp}, we have
\begin{eqnarray*}
g(P(\nabla^g_ey),z)=\omega(P(\nabla^g_ey),P(z))=-\omega(\nabla^g_ey,P^2(z))=-g(\nabla^g_ey,z)=0,
\end{eqnarray*}
which implies that
\begin{equation}\label{eq:ttt1}
g(\nabla^g_eP(y),z)=g(P(\nabla^g_ey),z),\quad \forall e\in\Gamma(E),y,z\in\Gamma(E^+).
\end{equation}
For all $e\in\Gamma(E),y\in\Gamma(E^+),\xi\in\Gamma(E^-)$, we have
\begin{eqnarray}\label{eq:ttt2}
g(\nabla^g_eP(y),\xi)=g(\nabla^g_ey,\xi)=-g(\nabla^g_ey,P(\xi))=g(P(\nabla^g_ey),\xi).
\end{eqnarray}
By \eqref{eq:ttt1} and \eqref{eq:ttt2}, we deduce that
\begin{eqnarray*}
 \nabla^g_{e}P(y) = P(\nabla^g_{e }y) ,\quad \forall e \in\Gamma(E),y\in\Gamma(E^+).
\end{eqnarray*}
Similarly, we have
\begin{eqnarray*}
\nabla^g_{e}P(\xi) = P(\nabla^g_{e }\xi) ,\quad \forall e \in\Gamma(E),\xi\in\Gamma(E^-).
\end{eqnarray*}
\eqref{eq:Levi-PC} follows immediately.\qed
\begin{cor}
For all $e\in\Gamma(E)$, we have $\nabla_e^g \Gamma(E^+)\subseteq \Gamma(E^+)$ and $\nabla_e^g \Gamma(E^-)\subseteq \Gamma(E^-).$
\end{cor}
\begin{pro}
Let $(E,[\cdot,\cdot]_E,\rho)$ be a Lie algebroid with a pseudo-Riemannian metric $g$. Let $P:E\longrightarrow E$ be a bundle map satisfying $P^2=\id $, \eqref{eq:Levi-PC} and \eqref{eq:Rie-PC}. Then there is a natural non-degenerate closed $2$-form defined by
\begin{equation}\label{eq:Rie-Symp}
\omega(e_1,e_2)=g(e_1,P(e_2)),\quad \forall e_1,e_2\in\Gamma(E).
\end{equation}
Furthermore, $(E,\omega,P)$ is a para-K\"{a}hler Lie algebroid.
\end{pro}

\pf It follows from  the same proof of Proposition \ref{pro:Levi-PC}. \qed

\subsection{Para-complex structures on pre-symplectic algebroids}

 \begin{defi}
 Let $(E,\star,\rho,(\cdot,\cdot)_-)$ be a pre-symplectic algebroid. A {\bf
para-complex  structure} on $E$ is a bundle map  $P:E\longrightarrow
E$ satisfying the algebraic conditions:
\begin{equation}\label{eq:compatiblity}
P^2={\id},\quad (P(x),P(y))_-=-(x,y)_-,\quad \forall x,y\in\Gamma(E),
\end{equation}
and the integrability  condition:
 \begin{equation}\label{eq:int-cond2}
P(e_1\star e_2)=P(e_1)\star e_2+e_1\star P(e_2)-P(P(e_1)\star P(e_2)),\quad\forall e_1,e_2\in\Gamma(E).
\end{equation}
We call a pre-symplectic algebroid $(E,\star,\rho,(\cdot,\cdot)_-)$ with a para-complex  structure $P$ a {\bf para-complex pre-symplectic algebroid}.
\end{defi}
Denote a para-complex pre-symplectic algebroid by $(E,\star,\rho,(\cdot,\cdot)_-,P)$.

Let us also denote
\begin{eqnarray*}
E^+_m:=\{x\in E_m\mid P_m(x)=x,~\forall m\in M\},\quad E^-_m:=\{\xi\in E_m\mid P_m(\xi)=-\xi,~\forall m\in M\}.
\end{eqnarray*}
By the algebraic condition \eqref{eq:compatiblity}, $E^+$ and $E^-$ are isotropic subbundles of vector bundle $E$ over $M$ and $E$ is the Whitney sum of $E^+$ and $E^-$, i.e. $E=E^+\oplus E^-$. By the integrability  condition \eqref{eq:int-cond2}, we have

\begin{lem}\label{lem:int-subalg}
 Let $P$ be a para-complex  structure on a pre-symplectic algebroid $(E,\star,\rho,(\cdot,\cdot)_-)$. Then   $E^+$ and $E^-$ are transversal Dirac structures of the pre-symplectic algebroid $(E,\star,\rho,(\cdot,\cdot)_-)$.
\end{lem}

 \begin{pro}\label{pro:PC-PK}
 Let $(E,\star,\rho,(\cdot,\cdot)_-,P)$ be a para-complex pre-symplectic algebroid. Then \\$(E,[\cdot,\cdot]_E,\rho,\omega=(\cdot,\cdot)_-,P)$ is a para-K\"{a}hler Lie algebroid.

 Conversely, if $(E,[\cdot,\cdot]_E,\rho,\omega,P)$ is a para-K\"{a}hler Lie algebroid, then $(E,\star,\rho,(\cdot,\cdot)_-=\omega,P)$ is a para-complex pre-symplectic algebroid.
 \end{pro}
\pf It is easy to see that the integrability  condition \eqref{eq:int-cond2} implies the integrability  condition \eqref{eq:int-cond1}. By   \eqref{eq:compatiblity}, we can deduce that \eqref{eq:Para-Comp} holds. Thus, $(E,[\cdot,\cdot]_E,\rho,\omega=(\cdot,\cdot)_-,P)$ is a para-K\"{a}hler Lie algebroid.

 Conversely, if $(E,[\cdot,\cdot]_E,\rho,\omega,P)$ is a para-K\"{a}hler Lie algebroid, by the fact that   $E^+$ and $E^-$ are two transversal Lagrangian subalgebroids, we can deduce that the integrability  condition \eqref{eq:int-cond2} holds. We omit details.  \qed

 \emptycomment{
 by Eq. \eqref{eq:Para-Comp}, we have
$\Gamma(E^+)$ is isotropic associated to $\omega$. Then for $x,y,z\in\Gamma(A^+)$, by Proposition \ref{Pro:int-equi}, we have
\begin{eqnarray*}
\omega(P(x\star y),z)&=&-\omega(x\star y,P(z))=-\omega(x\star y,z)\\
&=&-\rho(x)\omega(y,z)-\frac{1}{2}\rho(z)\2omega(x,y)+\omega(y,[x,z])=0,
\end{eqnarray*}
and similarly, we have $\omega(x\star y,z)=0$. Thus we have
$$\omega(P(x\star y),z)=\omega(x\star y,z),\quad x,y,z\in\Gamma(A^+).$$
For $x,y\in\Gamma(A^+),\xi\in\Gamma(A^-)$, we have
$$\omega(P(x\star y),\xi)=-\omega(x\star y,P(\xi))=\omega(x\star y,\xi).$$
Thus for all $x,y\in\Gamma(A^+)$, we have $P(x\star y)=x\star y$. Similarly, for all $\xi,\eta\in\Gamma(A^-)$, we have $P(\xi\star \eta)=-\xi\star \eta$.
For $x\in\Gamma(A^+),\xi\in\Gamma(A^-)$, we have
\begin{eqnarray*}
&&P(x\star \xi)-P(x)\star \xi-x\star P(\xi)+P(P(x)\star P(\xi))\\
&=&P(x\star \xi)-P(x\star \xi)=0
\end{eqnarray*}
Similarly, $P(\xi\star x)=P(\xi)\star x+\xi\star P(x)-P(P(\xi)\star P(x)).$
Thus integrability  condition \eqref{eq:int-cond2} holds.}

 \begin{pro}\label{pro:mul-Levi}
Let $(E,\star,\rho,(\cdot,\cdot)_-,P)$ be a para-complex pre-symplectic algebroid.  Then $(E,g)$ is a pseudo-Riemannian Lie algebroid, where $g$ is the pseudo-Riemannian metric defined by \eqref{eq:ParaC-Rie}, in which $\omega=(\cdot,\cdot)_-$.  Furthermore, for all $x,y\in\Gamma(E^+)$ and $\xi,\eta\in\Gamma(E^-)$, we have
$$\nabla^g_xy=x\star y,\quad \nabla^g_\xi\eta=\xi\star \eta,$$
where $\nabla^g$ is the Levi-Civita $E$-connection defined by \eqref{eq:nablag}.
\end{pro}
\pf By Proposition \ref{pro:PC-PK} and Proposition \ref{pro:Levi-PC}, $(E,g)$ is a pseudo-Riemannian Lie algebroid. For all $x,y\in\Gamma(E^+),~\xi\in\Gamma(E^-)$, by the fact that $E^+$ is isotropic and condition (ii) in Definition \ref{MTLA}, we have
\begin{eqnarray*}
(x\star y,\xi)_-&=&\rho(x)(y,\xi)_--(y,[x,\xi]_A)_-.
\end{eqnarray*}
By \eqref{eq:ParaC-Rie} and condition (i) in Definition \ref{defi:LCconnection}, we have
\begin{eqnarray*}
g(x\star y,\xi)&=&\rho(x)g(y,\xi) -g(y,[x,\xi]_A)=\rho(x)g(y,\xi) -g(y,\nabla^g_x\xi)=g(\nabla^g_xy,\xi).
\end{eqnarray*}
On the other hand, for $x,y,z\in\Gamma(E^+)$, we have $g(x\star y,z)=g(\nabla^g_xy,z)=0$.
 Thus, for all $x,y\in\Gamma(E^+)$, we have $\nabla^g_xy=x\star y$. Similarly, for all $\xi,\eta\in\Gamma(E^-)$, $\nabla^g_\xi\eta=\xi\star \eta$.\qed

\begin{rmk}
The above result means that
the pre-symplectic algebroid structure can characterize the restriction to Lagrangian subalgebroids of the Levi-Civita $E$-connection $\nabla^g$ in the corresponding pseudo-Riemannian Lie algebroid.
\end{rmk}

\begin{ex}{\rm
 Consider the pre-symplectic algebroid $(E=A\oplus A^*,\star,\rho,(\cdot,\cdot)_-)$ given in Proposition \ref{pro:standard}, where $\star$ is given by \eqref{eq:standardLWXbracket}, $(\cdot,\cdot)_-$ is given by \eqref{eq:defiomega} and $\rho=a_A$. The bundle map $P:E\longrightarrow E$ defined by
 \begin{equation}\label{eq:paracomplex}
 P(x+\xi)=x-\xi,\quad \forall x\in\Gamma(A),\xi\in\Gamma(A^*),
 \end{equation}
 is a para-complex  structure on $E$.
}
\end{ex}

\end{document}